\documentclass{amsart}
\usepackage{amsthm}
\usepackage{amssymb}
\newcounter{idee}

\swapnumbers
\newtheorem{theorem}[idee]{Theorem}

\newtheorem{corollary}[idee]{Corollary}
\newtheorem{lemma}[idee]{Lemma}

\theoremstyle{definition}

\newtheorem{comments}[idee]{Comment}

\theoremstyle{remark}
\newtheorem{remark}[idee]{Remark}
\newtheorem{notation}[idee]{Notation}
\newtheorem{notations}[idee]{Notations}

\begin{document}
\title{A note on idempotents in finite AW*-factors}
\author{Gabriel Nagy}
\address{Department of Mathematics, Kansas State university, Manhattan KS 66506, U.S.A.}
\email{nagy@math.ksu.edu}
\thanks{Partially supported by NSF grant DMS 9706858}
\keywords{AW*-algebra, quasi-trace, idempotent, projection, dimension function}
\subjclass{Primary 46L10; Secondary 46L30}
\begin{abstract}
We prove that the value of the quasi-trace on an idempotent
element in a AW*-factor of type $\text{II}_1$ is the same
as the dimension of its left (or right) support.\end{abstract}
\maketitle

It is a long standing open problem (due to Kaplansky) to prove that {\em
an AW*-factor of type $\text{II}_1$ is in fact a von Neumann algebra}.
A remarkable
answer, in the affirmative, was found by Haagerup (\cite{Ha}), who proved that
if an AW*-factor $A$ is generated by an exact C*-algebra,  then $A$ is indeed
a von Neumann algebra.

The main object, that was investigated in connection with Kaplansky's problem,
is the {\em quasi-trace}, whose construction we briefly recall below.

One starts with an AW*-factor of type $\text{II}_1$, say $A$. Denote by
$\mathcal{P}(A)$ the collection of projections in $A$, that is
$$
\mathcal{P}(A)=\{p\in A\,:\,p=p^*=p^2\}.$$
A key fact is then the existence of a (unique)
{\em dimension function} $D:\mathcal{P}(A)
\to [0,1]$ with the following properties:
\begin{itemize}
\item $D(p)=D(q)\Longleftrightarrow p\sim q$;
\item if $p\perp q$, then $D(p+q)=D(p)+D(q)$;
\item $D(1)=1$.
\end{itemize}
The symbol ``$\sim$'' denotes the Murray-von Neumann equivalence relation
($p\sim q\Leftrightarrow \exists\,x\in A$ with $p=x^*x$ and $q=xx^*$), while
``$\perp$'' denotes the orthogonality relation ($p\perp q\Leftrightarrow
pq=0$; this implies that $p+q$ is again a projection).

Once the dimension function is defined, it is extended to self-adjoint elements
with finite spectrum. More explicitly, if $a\in A$ is self-adjoint with
finite spectrum, then there are (real) numbers $\alpha_1,\dots,\alpha_n$
and pairwise orthogonal projections $p_1,\dots,p_n$, such that
$a=\sum_{k=1}^n\alpha_k p_k$. We then define $d(a)=\sum_{k=1}^n\alpha_k D(p_k)$.

For an arbitrary self-adjoint element $a\in A$,
one can approximate uniformly $a$ with a sequence $(a_n)_{n\geq 1}\in\{a\}''$
of elements with finite spectrum. (Here $\{a\}''$ stands for the AW*-subalgebra
generated by $a$ and $1$.) It turns out that the limit
$q(a)=\lim_{n\to\infty}d(a_n)$ is independent of the particular choice of
$(a_n)_{n\geq 1}$.

Finally, for an arbitrary element $x\in A$, one defines $Q(x)=q(\text{Re}\, x)+
iq(\text{Im}\, x)$, where $\text{Re}\, x=\frac 12(x+x^*)$ and $\text{Im}\, x=\frac 1{2i}(x-x^*)$.

The map $Q:A\to\mathbb{C}$, defined this way, is the unique one with the
properties:
\begin{itemize}
\item[(i)] $Q$ is {\em linear}, when restricted to {\em abelian\/} C*-subalgebras
of $A$;
\item[(ii)] $Q(x^*x)=Q(xx^*)\geq 0$, for all $x\in A$;
\item[(iii)] $Q(x)=Q(\text{Re}\, x)+
iQ(\text{Im}\, x)$,  for all $x\in A$;
\item[(iv)] $Q(1)=1$.
\end{itemize}
It is obvious that $Q\big|_{\mathcal{P}(A)}=D$. The map $Q$ is called {\em the
quasi-trace of} $A$.

It is well known that an AW*-factor of type $\text{II}_1$
is a von Neumann algebra, if and only if its quasi-trace is {\em linear}.
Haagerup's solution for Kaplansky's problem goes through the
proof of the linearity of the quasi-trace.

On the one hand, one can easily see that the linearity of the quasi-trace is
equivalent
to its scalar homogeneity (compare with (i) above):
\begin{equation}
Q(\alpha x)=\alpha Q(x)\text{, for all }x\in A,\,\,\alpha\in\mathbb{C}.
\tag{H}\end{equation}
Notice that (H) holds when either $\alpha\in\mathbb{R}$, or when $x$ is
{\em normal}.
On the other hand, it is again easy to note that the linearity
of the quasi-trace is equivalent to the similarity invariance property
\begin{equation}
Q(sxs^{-1})=Q(x)\text{, for all }x\in A,\,\,s\in GL(A).
\tag{S}\end{equation}
(Here $GL(A)$ denotes the group of invertible elements in $A$.)
Notice that (S) is true if $s$ is {\em unitary}.

The purpose of this note is to prove that {\em both\/} (H) {\em and\/} (S)
{\em hold, if $x\in A$ is an idempotent\/} (i.e. $x^2=x$).


\


\begin{notations} If $A$ is an AW*-algebra, for an element $x\in A$, we denote
by $\mathbf{L}(x)$ (resp. $\mathbf{R}(x)$) the {\em left\/} (resp. {\em right\/})
{\em support\/} of $x$. Recall that both $\mathbf{L}(x)$ and
$\mathbf{R}(x)$ are projections, and moreover we have
$$
\mathbf{L}(x)\sim\mathbf{R}(x)\text{, for all }x\in A.$$

It is known (see \cite{Ka}) that, for every $x\in A$, there exists a unique
partial isometry $v$ such that
\begin{itemize}
\item[(i)] $x= v(x^*x)^{1/2}$;
\item[(ii)] $vv^*=\mathbf{L}(x)$;
\item[(iii)] $v^*v=\mathbf{R}(x)$.
\end{itemize}
(This property is referred to as the {\em Polar Decomposition}.)
\end{notations}

\begin{remark} If $A$ is a finite AW*-algebra, then the group
$GL(A)$, of invertible elements, is dense in $A$ in the norm topology.

Indeed, on the one hand, since $\mathbf{L}(x)\sim\mathbf{R}(x)$, by the
finiteness assumption we also have
$1-\mathbf{L}(x)\sim 1-\mathbf{R}(x)$. In particular, there exists a
partial isometry $w\in A$ such that $1-\mathbf{L}(x)=ww^*$ and
$1-\mathbf{R}(x)=w^*w$. On the other hand,
if we take $x=v(x^*x)^{1/2}$ be the polar decomposition described
above, then we obviously have $w^*v=0$, so the element $u=v+w$ is
{\em unitary}, and we still have $x=u(x^*x)^{1/2}$. Then, for every
$\varepsilon>0$ the positive element $(x^*x)^{1/2}+\varepsilon 1$ is
invertible, and so is $u\{(x^*x)^{1/2}+\varepsilon 1\}$.
The result then follows from the obvious equality
$
\|x-u\{(x^*x)^{1/2}+\varepsilon 1\}\|=\varepsilon$.
\end{remark}

\begin{notation} For a C*-algebra $A$, and an integer $n\geq 2$, we denote
by $M_n(A)$ the C*-algebra of $n\times n$ matrices with coefficients in $A$.
\end{notation}

The key technical result in this paper is the following.

\begin{lemma} Suppose $A$ is a unital C*-algebra, and $x\in GL(A)$. Then there
exists a unitary element $U\in M_2(A)$, and elements $y,z\in A$, such that
\begin{equation}
U^*\begin{bmatrix} 2 & x \\ 0 & 0\end{bmatrix}U=\begin{bmatrix} 1 & y \\ z& 1
\end{bmatrix}\label{lema1eq}\end{equation}\end{lemma}

\begin{proof} Consider the function defined by
$$
f(t)=\frac{(1-t)^2}t,\,\,\,0<t<1.$$
It is obvious that $f:(0,1)\to (0,\infty)$ is a homeomorphism. Since $x$ is invertible,
the spectrum of $xx^*$ is contained in $(0,\infty)$, so by functional calculus there
exists an invertible positive element $w\in A$, with $\|w\|<1$, such that
$f(w)=xx^*$, which means
\begin{equation}
xx^*=(1-w)^2w^{-1}.\label{lemma1-w}\end{equation}

Define the elements
\begin{align*}
a&=(1+w)^{-1/2}w^{1/2}, &b&=(1+w)^{-1/2},\\
c&=x^{-1}(1+w)^{-1/2}w^{-1/2}(1-w), &d&=-x^{-1}(1+w)^{-1/2}(1-w).
\end{align*}
First, we have
\begin{equation}
aa^*+bb^*=(1+w)^{-1/2}[w+1](1+w)^{-1/2}=1,\label{aa*}\end{equation}
and using \eqref{lemma1-w} we also have
\begin{align}
cc^*+dd^*&=x^{-1}(1+w)^{-1/2}(1-w)[w^{-1}+1](1-w)(1+w)^{-1/2}(x^*)^{-1}=\label{cc*}\\
&=x^{-1}(1-w)^2w^{-1}(x^*)^{-1}=x^{-1}(xx^*)(x^*)^{-1}=1;\notag\end{align}
Secondly, since by taking inverses, \eqref{lemma1-w} yields
\begin{equation}
(x^*)^{-1}x^{-1}=w(1-w)^{-2},\label{lemma1-winv}\end{equation}
so we also get
\begin{align}
a^*a+c^*c&=w^{1/2}(1+w)^{-1/2}[1+(1-w)(x^*)^{-1}x^{-1}(1-w)](1+w)^{-1/2}w^{1/2}=
\notag\\
&=w^{1/2}(1+w)^{-1/2}[1+w^{-1}](1+w)^{-1/2}w^{1/2}=1;
\label{a*a}\\
b^*b+d^*d&=(1+w)^{-1/2}[1+(1-w)(x^*)^{-1}x^{-1}(1-w)](1+w)^{-1/2}=
\notag\\
&=(1+w)^{-1/2}[1+w](1+w)^{-1/2}1.
\label{b*b}
\end{align}
Finally, we notice that
\begin{equation}
ac^*+bd^*=(1+w)^{-1/2}[w^{1/2}w^{-1/2}-1](1-w)(1+w)^{-1/2}(x^*)^{-1}=0,\label{ac*}
\end{equation}
and using \eqref{lemma1-winv} we also have
\begin{align}
a^*b+c^*d&=(1+w)^{-1/2}[w^{1/2}-w^{-1/2}(1-w)(x^*)^{-1}x(1-w)](1+w)^{-1/2}=\label{a*b}\\
&=(1+w)^{-1/2}[w^{1/2}-w^{-1/2}w](1+w)^{-1/2}=0.\notag
\end{align}
If we define the matrix $U=\begin{bmatrix} a& b \\ c& d\end{bmatrix}$, then
\eqref{aa*}, \eqref{cc*}, and \eqref{ac*} give $UU^*=I$, while
\eqref{a*a}, \eqref{b*b}, and \eqref{a*b} give $U^*U=I$, so $U$ is indeed unitary.
(Here $I$ denotes the unit in $M_2(A)$.)

Let us observe now that
\begin{align}
2a+xc&=2(1+w)^{-1/2}w^{1/2}+(1+w)^{1/2}(1-w)w^{-1/2}=
\label{11}\\
&=(1+w)^{-1/2}(w^{1/2}+w^{-1/2})=a+bw^{-1/2};
\notag\\
2b+xd&=2(1+w)^{-1/2}-(1+w)^{-1/2}(1-w)=
\label{12}\\
&=(1+w)^{-1/2}(w+1)=aw^{1/2}+b;
\notag\\
c+dw^{-1/2}&=x^{-1}(1+w)^{-1/2}[(1-w)w^{-1/2}-(1-w)w^{-1/2}]=0;\label{21}\\
cw^{1/2}+d&=(c+dw^{-1/2})w^{1/2}=0.\label{22}
\end{align}
These equalities prove exactly that
$$
\begin{bmatrix} 2 & x\\ 0&0\end{bmatrix}\cdot
\begin{bmatrix} a & b\\ c&d\end{bmatrix}=
\begin{bmatrix} 1 & w^{1/2}\\ w^{-1/2} & 1\end{bmatrix}\cdot
\begin{bmatrix} a&b \\ c&d\end{bmatrix}.$$
\end{proof}

The next result is a particular case of the main result

\begin{lemma}
Let $A$ be an AW*-factor of type $\text{\rm II}_1$, and let $e\in
A$ be an idempo\-tent with $D(\mathbf{L}(e))=\frac 12$. Then, for any $\lambda\in\mathbb{C}$,
we have
\begin{equation}
Q(\lambda e)=\lambda /2.\end{equation}
\end{lemma}

\begin{proof}
Denote, for simplicity, the projection $\mathbf{L}(x)$ by $p$. The assumption is that
$p\sim 1-p$. Then we have a $*$-isomorphism $\Phi:A\to M_2(pAp)$, such that
$\Phi(p)=\begin{bmatrix} 1&0\\ 0&0\end{bmatrix}$. It is obvious that, since
$e=pe$, there exists an element $x\in pAp$ such that
$\Phi(e)=\begin{bmatrix} 1&0\\ x&0\end{bmatrix}$.
By Remark 2, we can find a sequence $(x_n)_{n\geq 1}$ of invertible elements
in $pAp$ (the unit in $pAp$ is $p$), with $\lim_{n\to\infty}\|x_n-x\|=0$. Define
the sequence
$$E_n=\begin{bmatrix} 1&x_n\\ 0&0\end{bmatrix}\in M_2(pAp),\,\,\,n\geq1.$$
By Lemma 4, one can find two sequences $(y_n)_{n\geq 1}$ and $(z_n)_{n\geq 1}$
in $pAp$, and a sequence of unitaries $(U_n)_{n\geq 1}\subset M_2(pAp)$, such that
\begin{equation}
E_n=U_n\begin{bmatrix} \tfrac 12& y_n\\ z_n&\tfrac 12\end{bmatrix}U_n^*
\text{, for all }n\geq 1.\label{lemma2-1}\end{equation}
Define $e_n=\Phi^{-1}(E_n)$, $u_n=\Phi^{-1}(U_n)$, and
$a_n=u_n^*e_nu_n-\frac 12 1$, $n\geq 1$, so that we have
\begin{equation}
\Phi(a_n)=\begin{bmatrix} 0&y_n\\ z_n&0\end{bmatrix}
\text{, for all }n\geq 1.\label{lemma2-2}\end{equation}
Fix now a complex number $\lambda$. On the one hand, we have
$$
\lambda e_n=u_n(\tfrac{\lambda}2+a_n)u_n^*=\tfrac{\lambda}2 1+u_na_nu_n^*
\text{, for all }n\geq 1.$$
This gives, for every $n\geq 1$, the equalities
\begin{align}
\text{Re}(\lambda e_n)&=(\text{Re}\tfrac{\lambda}2)1+\text{Re}(\lambda u_na_nu_n^*);
\label{lemma2-re}\\
\text{Im}(\lambda e_n)&=(\text{Im}\tfrac{\lambda}2)1+\text{Im}(\lambda u_na_nu_n^*).
\label{lemma2-im}
\end{align}
Notice however that, using the {\em unitary invariance\/} (property (S) for
$s$ unitary), together with \eqref{lemma2-re} and \eqref{lemma2-im}, gives,
for every $n\geq 1$, the equalities
\begin{align}
Q(\lambda e_n)&=Q\big((\text{Re}\tfrac{\lambda}2)1+\text{Re}(\lambda u_na_nu_n^*)\big)+
iQ\big((\text{Im}\tfrac{\lambda}2)1+\text{Im}(\lambda u_na_nu_n^*)\big)=\notag\\
&=\text{Re}\tfrac{\lambda}2+Q\big(\text{Re}(\lambda u_na_nu_n^*)\big)+
i\text{Im}\tfrac{\lambda}2+iQ\big(\text{Im}(\lambda u_na_nu_n^*)\big)=\label{lemma2-3}\\
&=\tfrac{\lambda}2+Q((\lambda u_n a_nu_n^*)=
\tfrac{\lambda}2+Q\big(u_n(\lambda a_n)u_n^*\big)=
\tfrac{\lambda}2+Q\big(\lambda a_n).\notag
\end{align}

On the other hand, if we define $v=1-2p$ (which is obviously a unitary in $A$), then
for all $n\geq 1$ we have
$$
\Phi\big(v(\lambda a_n)v^*\big)=
\begin{bmatrix} 1 & 0 \\ 0 &-1\end{bmatrix}\cdot
\begin{bmatrix}0 &\lambda y_n\\ \lambda z_n &0\end{bmatrix}
\cdot\begin{bmatrix} 1 & 0 \\ 0 &-1\end{bmatrix}=
\begin{bmatrix} 0 & -\lambda y_n \\ -\lambda z_n &0\end{bmatrix}=
\Phi(-\lambda a_n)
.$$
This gives
\begin{equation}
v(\lambda a_n)v^*=-\lambda a_n\text{, for all }n\geq 1.
\label{lemma2-v}\end{equation}
Using the unitary invariance, combined with {\em real scalar homogeneity\/}
(property (H) with $\alpha\in\mathbb{R}$), the equality \eqref{lemma2-v} gives
$$
Q(\lambda a_n)=Q(v(\lambda a_n)v^*)=Q(-\lambda a_n)=-Q(\lambda a_n),$$
which forces
\begin{equation}
Q(\lambda a_n)=0\text{, for all }n\geq 1.\label{lemma2-qa=0}\end{equation}
Combining \eqref{lemma2-qa=0} with \eqref{lemma2-3} gives
\begin{equation}
Q(\lambda e_n)=\lambda /2\text{, for all }n\geq 1.
\label{lemma2-qen}\end{equation}
It is obvious that, by construction, we have $\lim_{n\to\infty}\|E_n-\Phi(e)\|=0$,
which means that $\lim_{n\to\infty}\|e_n-e\|=0$. Using the norm continuity of the
quasi-trace (see \cite{BH}), combined with \eqref{lemma2-qen} gives the desired
result.
\end{proof}

We are now ready to prove the main result.

\begin{theorem} Let $A$ be an AW*-factor of type $\text{\rm II}_1$, and let
$e\in A$ be an idempotent element. Then, for any $\lambda\in\mathbb{C}$, one has
the equality
\begin{equation}
Q(\lambda e)=\lambda D\big(\mathbf{L}(e)\big).\label{thm-eq}\end{equation}
\end{theorem}

\begin{proof} The proof will be carried on in two steps.

{\sc Particular Case:} {\em Assume $D\big(\mathbf{L}(e)\big)\leq\frac 12$.}

Denote, for simplicity $\mathbf{L}(e)$ by $p$, and $\mathbf{R}(e)$ by $q$. By the
{\em Paralellogram Law\/} (see \cite{Ka}) we have
$$
p\vee q -p\sim q-p\wedge q.$$
Since $p\sim q$, we get
$$D(p\vee q-p)=D(q-p\wedge q)=D(q)-D(p\wedge q)\leq D(q)=D(p).$$
Using the intermediate value property for $D$, we can find a projection
$r\leq 1-p\vee q$ such that $D(r)+D(p\vee q-p)=D(p)$. Put
$q_0=r+p\vee q-p$. We have $p\perp q_0$, and $p+q_0\geq p\vee q$.
Let us work in the AW*-algebra $A_0=(p+q_0)A(p+q_0)$. Obviously $A_0$ is
again an AW*-factor of type $\text{II}_1$, so it carries its quasi-trace
$Q_0$. By the uniqueness of the quasi-trace, it is obvious that
\begin{equation}
Q_0(x)=\frac{Q(x)}{D(p+q_0)}\text{, for all }x\in A_0.
\label{thm-q0}\end{equation}
Notice that $e\in (p\vee q)A(p\vee q)$, so in particular $e$ belongs
to $A_0$. In $A_0$, we have $\mathbf{L}(e)=p$, and $p\sim q_0=1_{A_0}-p$, which
means that $D_0(p)=\frac 12$. (Here $D_0$ denotes the dimension function in
$A_0$). By Lemma 5, we get
$$Q_0(\lambda e)=\lambda/2,$$
which combined with \eqref{thm-q0} yields
$$
\frac{\lambda}2=\frac{Q(\lambda e)}{D(p+q_0)}=
\frac{Q(\lambda e)}{2D(p)},$$
which obviously proves \eqref{thm-eq}.

{\sc General Case}. One knows (see \cite{Be}) that $M_2(A)$ is also an
AW*-factor of type $\text{II}_1$. Moreover, if we denote by $Q^{(2)}$ the
quasi-trace of $M_2(A)$, we have (as above)
\begin{equation}
Q(x)=2 Q^{(2)}\big(\begin{bmatrix} x& 0 \\ 0 & 0\end{bmatrix}\big)
\text{, for all }x\in A.\label{thm-q2}\end{equation}
Define the idempotent
$E=\begin{bmatrix} e& 0 \\ 0 & 0\end{bmatrix}\in M_2(A)$, and the
projection $P=\begin{bmatrix} p& 0 \\ 0 & 0\end{bmatrix}\in M_2(A)$.
Using \eqref{thm-q2} we have $Q(e)=2Q^{(2)}(E)$. If we denote by $D^{(2)}$ the
dimension function of $M_2(A)$, then by \eqref{thm-q2} we also have
$D(p)=2D^{(2)}(P)$, which gives $D^{(2)}(P)\leq\frac 12$. Using the obvious
equality $\mathbf{L}(E)=P$, by the particular case above applied to $M_2(A)$,
we get
$$
Q(\lambda e)=2Q^{(2)}(\lambda E)=2\lambda D^{(2)}(P)=\lambda D(p).
$$
\end{proof}

\begin{corollary} Let $A$ be an AW*-factor of type $\text{\rm II}_1$, and let
$e\in A$ be an idempotent. Then, for every $s\in GL(A)$, we have
\begin{equation}
Q(ses^{-1})=Q(e).
\label{cor1-eq}\end{equation}
\end{corollary}

\begin{proof} Put $p=\mathbf{L}(e)$. If we define $t=1-e(1-p)$, then
$t$ is invertible, and $e=tpt^{-1}$. This computation shows that it is
enough to prove \eqref{cor1-eq} in the case when $e=e^*$.
Let $q=\mathbf{L}(ses^{-1})$. Arguing as above, there exists some $x\in GL(A)$
such that $ses^{-1}=xqx^{-1}$. So if we put $y=x^{-1}s$, we have
$yey^{-1}=q$, with both $e$ and $q$ self-adjoint idempotents. Since this
obviously forces $e\sim q$, we get $D(e)=D(q)=Q(ses^{-1}$
\end{proof}

\begin{corollary} Let $A$ be an AW*-factor of type $\text{\rm II}_1$ and let
$e_1,e_2\in A$ be idempotents, such that $e_1e_2=e_2e_1=0$. {\rm (This implies that
$e_1+e_2$ is again an idempotent.)} Then
$$
Q(e_1+e_2)=Q(e_1)+Q(e_2).$$
\end{corollary}

\begin{proof} Let $p=\mathbf{L}(e_1+e_2)$. As in the proof of the preceding Corollary,
there exists an element $s\in GL(A)$, such that $p=s(e_1+e_2)s^{-1}$.
Put $f_k=se_ks^{-1}$, $k=1,2$, so that $f_1+f_2=p$. But now, we also have
$pf_k=f_kp=f_k$, $k=1,2$, which means in particular that $f_1,f_2\in pAp$.
So, if we work in the AW*-factor (again of type $\text{II}_1$) $A_0=pAp$,
we will have $f_1+f_2=1_{A_0}$. On the one hand,
using the notations from the proof of Theorem 6,
we have
$$Q_0(f_1)+Q_0(f_2)=1.$$
On the other hand, we have
$$
Q_0(f_k)=\frac{Q(f_k)}{D(p)},\,\,\,k=1,2,$$
so we get
$$D(p)=Q(f_1)+Q(f_2).$$
Finally, using the preceding Corollary, we get
$$
Q(e_1+e_2)=D(p)=Q(f_1)+Q(f_2)=Q(se_1s^{-1})+Q(se_2s^{-1})
=Q(e_1)+Q(e_2).$$
\end{proof}

\begin{comments} A.
Theorem 6 can be analyzed from a different point of view, as follows. In principle,
one can extend the dimension function $D$ to the collection of all idempotents,
by $\tilde{D}(e)=D\big(\mathbf{L}(e)\big)$. One can easily prove that
this extended dimension function will have the same properties as the usual
dimension function (when Murray-von Neumann equivalence is extended to idempotents).
The point of Theorem 6 is then the fact that $\tilde{D}=Q$.

B. It is interesting to note that the linearity of the quasi-trace is equivalent
to the following condition:
\begin{itemize}
\item[($*$)] {\em For any family of idempotents $e_1,\dots,e_n\in A$, such that
$e_je_k=0$ for all $j,k\in\{1,\dots,n\}$ with $j\neq k$, and any family of
numbers $\alpha_1,\dots,\alpha_n\in\mathbb{R}$, one has:}
$$Q(\alpha_1e_1+\dots +\alpha_n e_n)=\alpha_1Q(e_1)+\dots +\alpha_nQ(e_n).$$
\end{itemize}
This will be discussed in a future paper.
\end{comments}

\end{document}